\documentclass{IEEEtran}
\usepackage{cite}
\usepackage{amsmath,amssymb,amsfonts}
\usepackage{algorithmic}
\usepackage{graphicx}
\usepackage{color}
\usepackage{textcomp}
\def\BibTeX{{\rm B\kern-.05em{\sc i\kern-.025em b}\kern-.08em
    T\kern-.1667em\lower.7ex\hbox{E}\kern-.125emX}}
\begin{document}
\title{Global Complex Roots and Poles Finding Algorithm Based on Phase Analysis for Propagation and Radiation Problems}
\author{Piotr Kowalczyk

\thanks{This work has been submitted to the IEEE for possible publication. Copyright may be transferred without notice, after which this version may no longer be accessible.}

\thanks{This work was supported under funding for Statutory Activities for the Faculty of Electronics, Telecommunication and Informatics,
Gdansk University of Technology.}
\thanks{Piotr Kowalczyk is with Gdansk University of Technology, Faculty of Electronics, Telecommunications and Informatics, Narutowicza 11/12, 80233 Gdansk, Poland (e-mail: pio.kow@gmail.com, piotr.kowalczyk@pg.gda.pl). }
}

\maketitle

\begin{abstract}
A flexible and effective algorithm for complex roots and poles finding is presented. A wide class of analytic functions can be analyzed, and any arbitrarily shaped search region can be considered. The method is very simple and intuitive. It is based on sampling a function at the nodes of a regular mesh, and on the analysis of the function phase. As a result, a set of candidate regions is created and then the roots/poles are verified using a discretized Cauchy's argument principle. The accuracy of the results can be improved by the application of a self-adaptive mesh. The effectiveness of the presented technique is supported by numerical tests involving different types of structures, where electromagnetic waves are guided and radiated. The results are verified, and the computational efficiency of the method is examined. 
\end{abstract}


\section{Introduction}
\label{sec:introduction}

\IEEEPARstart{M}{any} propagation and radiation problems are formulated in a complex domain. One of the most common parameter representing electromagnetic wave propagation, as guiding, radiation or losses, is in general a complex number. Similarly, resonant frequencies of resonators or radiators are complex values. In many cases the evaluation of these parameters boils down to finding a complex root of a more or less complicated function. Even in the case of a simple technique, as mode matching or field matching \cite{warecka18}, the roots cannot be found analytically. For more sophisticated structures and methods the function can be expressed in the numerical form (spectral domain approach, discrete methods, nonlinear matrix eigenvalue problem) \cite{lech14,lech18}. Therefore, effective and efficient root finding algorithms became a necessary tool in the electromagnetic waves analysis. 

Root finding is one of the oldest and most common numerical problems. For a real function of a real variable the problem can be solved using many different techniques. Moreover, the results can be verified simply by checking the function sign changes at the ends of any sufficiently small interval containing a single root. However, even in this case, finding all the roots in a fixed region can be difficult. Finding the roots of a complex valued function of a complex variable is more complicated. Although, there are many complex root finding techniques, they usually consider only a special class of functions or a restricted region of analysis. 

The most common schemes, such as Newton's method \cite{Abramowitz72} or Muller's method \cite{Press92} are useful if the initial value of the root is already known. The same applies to algorithms tracking the root in a function with an extra parameter \cite{Michalski11,Kowalczyk17AP}. Global root finding algorithms are very efficient for simple polynomial functions \cite{Pinkert76,Schonhage82}, so a number of procedures based on polynomial approximation has been proposed \cite{Long98,Changying10,KRAVANJA2000}. However, the zeros of the considered function may bear little or no relation to its polynomial approximation \cite{Delves67}. An extreme example is $f(z)=\exp(z)$, which has no roots, whereas its any finite polynomial approximation e.g. $f(z)\approx \sum_{n=0}^{N} z^n/n!$ has $N$ zeros. Moreover, the roots can be extremely sensitive to perturbations in the coefficients of the higher order polynomial \cite{Wilkinson1994}. Therefore, the results obtained from such approximation should be further verified. Furthermore, the accuracy of the obtained zeros cannot be simply determined and controlled, so they can only be used as starting values for an extra iterative process. However, the results of such iterative techniques can be unreliable, especially if many roots/poles are located in a small region. In such cases, there is no guarantee that the process will converge to the specific initially evaluated root. Therefore, some of the results can be omitted (if their initial values are not sufficiently close to the roots). The control of such processes can be a fairly difficult task. 

Moreover, the polynomial approximation is ineffective for the function containing singularities and branch cuts in the analyzed region (the same limitation applies to some mesh based methods \cite{Wan11,Meylan2003}). Instead, a rational approximation can be applied. An interesting brief review of the polynomial and rational approximations used for roots/poles determination (based on roots of unity disk) can be found in \cite{Austin2014}. The approach involving rational approximation seems to be much more flexible and accurate, however, it is also more fragile. The improvement is often at the expanse of generating spurious poles-zeros pairs (Froissart doublets). In some cases the problem can reduced by proper regularization \cite{Gonnet2011}, but still the roots/poles obtained from such approximation should be verified and their accuracy cannot be simply controlled. 

The previously mentioned methods can be very efficient, especially for simple functions, however, in many practical applications, difficult and complex routines are implemented (e.g., those based on a genetic algorithm \cite{Tian09,Ariyaratne16} or on knowledge of the function singularities \cite{Chen2017}). Recently, two global algorithms has been proposed \cite{Kowalczyk15,ZOUROS2018}, which are general and flexible. Although they are quite effective, their efficiency and reliability can be significantly improved.

In this article, a simple global complex roots and poles finding algorithm is presented. The technique can be applied for very wide class of analytic functions (including those containing singularities or even branch cuts). An arbitrarily shaped search region can be considered, so an extra numerical error (corresponding to scaling or mapping of the function) can be avoided. 

In the first step, the function is sampled using a regular triangular mesh. The idea of domain triangulation for finding the zeros of the function is not new. Its origins are rooted in multidimensional bisection \cite{Eiger84} and it is also used in \cite{Kowalczyk15}. However, in the presented technique the function phase in the nodes is analyzed, rather than the simple sign change. From this analysis "candidate edges" are detected. Next, all the triangles attached to the candidate edges are surrounded by close contours determining the "candidate regions". For these contours a discretized form of Cauchy's Argument Principle (CAP) is applied, in order to verify the existence of roots or poles in the candidate regions. The Discretized Cauchy's Argument Principle (DCAP) does not require the derivative of the function and integration over the contour, as it is presented in \cite{Kravanja1999,Henrici88} and \cite{GILLAN2006}. In the proposed approach a minimal number of the function samples is utilized for DCAP (sometimes only four) and the contour shape is determined by the mesh geometry.

To improve the accuracy of the results any local (iterative) root finding scheme can be applied. However, as previously indicated, such methods may be unreliable and some of the roots can be missed, if the initial value is not sufficiently precise. In the presented approach, a simple self-adaptive mesh refinement (inside the previously determined candidate regions) is applied. This approach has slightly worse convergence than three-point iterative algorithms (e.g. Muller's technique), but the results are much more reliable - if the root/pole is located inside the candidate region, it cannot be lost in the sequential iterations.

The proposed technique consists of two stages: the preliminary estimation and the self-adaptive mesh refinement. The latter stage can be skipped, if the required accuracy is obtained in the former stage. 

In order to support the validity of the presented technique several numerical tests, involving different types of functions, are performed. The results are verified using other global techniques \cite{Kowalczyk15,Austin2014,ZOUROS2018} and the computational effectiveness and efficiencies of the methods are compared. It is shown that the proposed algorithm can be up to three orders of magnitude faster and requires significantly smaller number of function evaluations.

The examples presented in this paper are focused on microwave and optical applications, however the algorithm is not limited to computational electrodynamics. Similar problems occur in acoustics \cite{Jensen2011}, control theory \cite{Popov62}, quantum mechanics \cite{Fernandez2001} and many other fields. 

\section{Algorithm}

Let us denote the considered analytic function by $f(z)$ and the search region by
$\Omega\subset \mathbb{C}$. The aim is to find all the zeros and poles of the function in this region. 

The proposed algorithm can be divided into two separate stages: preliminary
estimation and final refinement. In the preliminary estimation process the roots and poles are initially determined by sampling the function at the nodes of a triangular regular mesh and by using DCAP. In the second stage a self-adaptive mesh refinement is applied to obtain the required accuracy.

\subsection{Preliminary Estimation}
To increase the readability and clarity of the description of this stage, it is divided into five steps.

\subsubsection{Mesh}
In the first step, region $\Omega$ is covered with a regular triangular mesh (e.g. using Delaunay triangulation) of $N$ nodes and $P$ edges. A honeycomb arrangement (equilateral triangles) of the nodes $z_n\in \Omega$, results in the highest efficiency of the algorithm. However, any other configuration is also possible, provided the longest edge length is smaller than the assumed resolution $\Delta r$ (this length is discussed in more detail in section \ref{sec:limits}). 

\subsubsection{Function Evaluation}
In this step, the function is evaluated at all the nodes $f_n=f(z_n)$ (this part of the algorithm can be simply parallelized, which can significantly improve the efficiency of the process for large problems). In this method the function value is not as important as the quadrant in which it lies, and only the quadrant 
\begin{equation}
\label{eqn:quadrant}
Q_n=
\left\{
\begin{array}{ll}
1,&0\le \arg f(z_n)<\pi/2\\
2,&\pi/2\le \arg f(z_n)<\pi\\
3,&\pi \le \arg f(z_n)<3\pi/2\\
4,&3\pi/2 \le \arg f(z_n)<2\pi
\end{array}
\right.
\end{equation}
associated with the node will be taken into account in the subsequent part of the algorithm. 

\subsubsection{Candidate Edges}
Next, the phase change along each of the edges is analyzed. For this purpose, an extra parameter representing the quadrant difference along the edge can be introduced 
\begin{equation}
\label{eqn:Ep}
\Delta Q_p=Q_{n_{p2}}-Q_{n_{p1}},\qquad \Delta Q_p\in\{-2,-1,0,1,2\}
\end{equation}
where $n_{p1}$ and $n_{p2}$ are nodes attached to edge $p$. 

\begin{figure}
\centerline{\includegraphics[width=1.0\columnwidth]{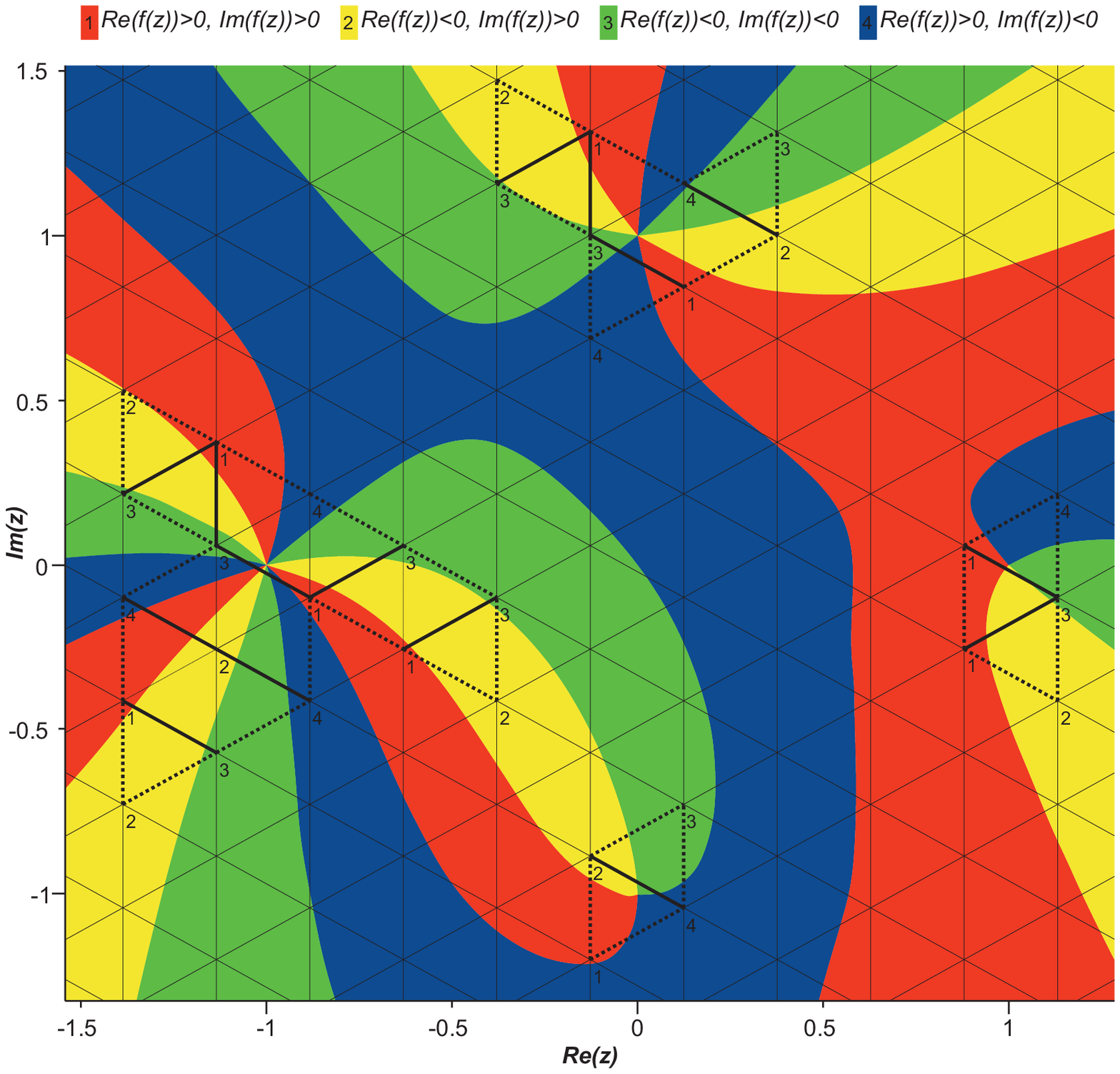}}
\caption{The preliminary estimation algorithm applied for function $f(z)=(z-1)(z-i)^2(z+1)^3/(z+i)$. The numbers (colors): $1$ (red), $2$ (yellow), $3$ (green) and $4$ (blue) represent the quadrants in which the function values lie. The candidate edges are marked by thick black lines. The black dotted lines represent the boundaries of the candidate regions.} \label{fig:E1}
\end{figure}

The main idea of this stage is based on the simple fact that any root or pole is located at the point where the regions described by four different quadrants meet - as it is shown in Figure~\ref{fig:E1} (to clarify this idea a phase portrait of the function is placed in the background \cite{Wegert2012}). Since any triangulation of the four nodes located in the four different quadrants requires at least one edge of $|\Delta Q_p|= 2$, then all such edges should be considered as a potential vicinity of the root or pole. All such candidate edges are collected in a single set 
 $\mathcal{E}=\{p: |\Delta Q_p|= 2\}$. 

\subsubsection{Candidate Regions}
All the triangles attached to the candidate edges $\mathcal{E}$ can also be collected in a single set of candidate triangles. From all the edges attached to these candidate triangles it is easy to find those that occur only once, and to collect them in a set $\mathcal{C}$, representing the boundary of the candidate regions (the inside edges are attached to two candidate triangles). The boundary of the candidate region must be constructed from the edges of   $|\Delta Q_p|<2$ only, as it is explained in the next paragraph \ref{VDCAP} (see condition (\ref{eqn:cond})).

Than, the set $\mathcal{C}$ can be divided into subsets $\mathcal{C}^{(k)}$, where $\mathcal{C}^{(k)}$ creates close contour surrounding $k$-th candidate region (in Figure \ref{fig:E1} there are four candidate regions). From the implementation point of view, such an operation is very simple. Starting from any edge from the set $\mathcal{C}$, one can construct the boundary of the region by finding the next edge connected to the previous one. If there is no connected edge in the set, then the edge should close the contour and the construction of the next candidate region can be started.

\subsubsection{Verification with Discretized Cauchy's Argument Principle}\label{VDCAP}
In a complex domain, to confirm the existence of a root or a pole in a fixed region, CAP is usually applied \cite{Brown2009}. According to this principle, the integral
\begin{equation}
\label{eqn:cap} q = \frac{1}{2 \pi i}\oint\limits_{C} \frac{f'(z)}{f(z)} dz
\end{equation}
represents the sum of all zeros counted with their multiplicities, minus the sum of all poles
counted with their multiplicities. If the region contains only a single candidate point, the parameter $q$ can be: a positive integer (root of order $q$), a negative integer (pole of order $-q$) or zero (regular point).

In practice, integral (\ref{eqn:cap}) represents a total change in the argument of the function $f(z)$ over a closed contour $C$ and there is no need to calculate this integral directly. The parameter $q$ can be evaluated from DCAP  - by sampling the function along the contour $C$ \cite{Kravanja1999,Henrici88} 
\begin{equation}
\label{eqn:dcap} q = \frac{1}{2\pi} \sum_{p=1}^P \textrm{arg} \frac{f(z_{p+1})}{f(z_{p})}.
\end{equation}
The points $z_{1},z_{2},...,z_{P}$ (and $z_{P+1}=z_{1}$) are obtained from discretization of the contour $C$ and the increase of the argument of $f(z)$ along the segment $C_p$ ($C=\bigcup_{p=1}^P C_p$) from $z_{p}$ to $z_{p+1}$ satisfies the condition
\begin{equation}
\label{eqn:cond} \left| [  \textrm{arg} f(z) ]_{z\in(z_{p},z_{p+1})} \right| \le \pi.
\end{equation}
As it is shown in \cite{Henrici88,Ying1988}, the condition (\ref{eqn:cond}) may not be easy to verify. However, in the presented approach the verification contour $C$ is defined by the boundary of the analyzed candidate region $\mathcal{C}^{(k)}$. For all the edges in $\mathcal{C}^{(k)}$ the phase change is $|\Delta Q_p|\le 1$ and the condition (\ref{eqn:cond}) is fulfilled.

An example of DCAP (single root in $z^{(1)}=1$, double root in $z^{(2)}=i$, triple root in $z^{(3)}=-1$ and singularity in $z^{(4)}=-i$) is presented in Figure~\ref{fig:E1}. For each of the four candidate regions, the function argument varies along the contour taking the values from the four quadrants. Since the quadrant difference along a single edge is $|\Delta Q_p|\le 1$, the discretization of the boundary is sufficient to evaluate the total phase change over the region boundary. By summing all the increases in the quadrants along the contour in the counterclockwise direction, one obtains the values $4$, $8$, $12$ and $-4$ for regions containing $z^{(1)}$, $z^{(2)}$, $z^{(3)}$ and $z^{(4)}$, respectively. Since the increase in the quadrant numbers along the edge represents the changes in the function argument of $\pi/2$, the parameter $q$ is equal to $1$, $2$ , $3$ and $-1$, respectively (single root, double root, triple root and singularity):
\begin{equation}
\label{eqn:cap_simp} q = \frac{1}{4}\sum_{p=1}^{P} \Delta Q_p.
\end{equation}

In general, at least $P=4q$ nodes is required to verify a single root or a pole of the $q$-th order.

It is worth noting that the change in quadrants along the candidate edges $|\Delta Q_p|=2$ is not unambiguous; it is impossible to determine whether the phase increases or decreases by two quadrants (condition (\ref{eqn:cond}) is not satisfied). 

In some cases, the regions cannot be unambiguously determined because the boundary of the candidate region cannot be closed (the candidate edge is located at the boundary of the domain $\Omega$). To solve this problem, the domain $\Omega$ should be extended or a denser initial mesh should be used.

\subsection{Mesh Refinement}
In order to improve the accuracy of the root location, a self-adaptive mesh is applied. This approach prevents an improper convergence of the algorithm - none of the initially found roots or poles can be missed even if they are not exactly inside the candidate region. In other techniques (such as Newton's or Muller's) a process that started with a given initial point can converge to a different root/pole, especially if the roots/poles are located in the immediate vicinity of each other.

\begin{figure}
\centerline{\includegraphics[width=1.0\columnwidth]{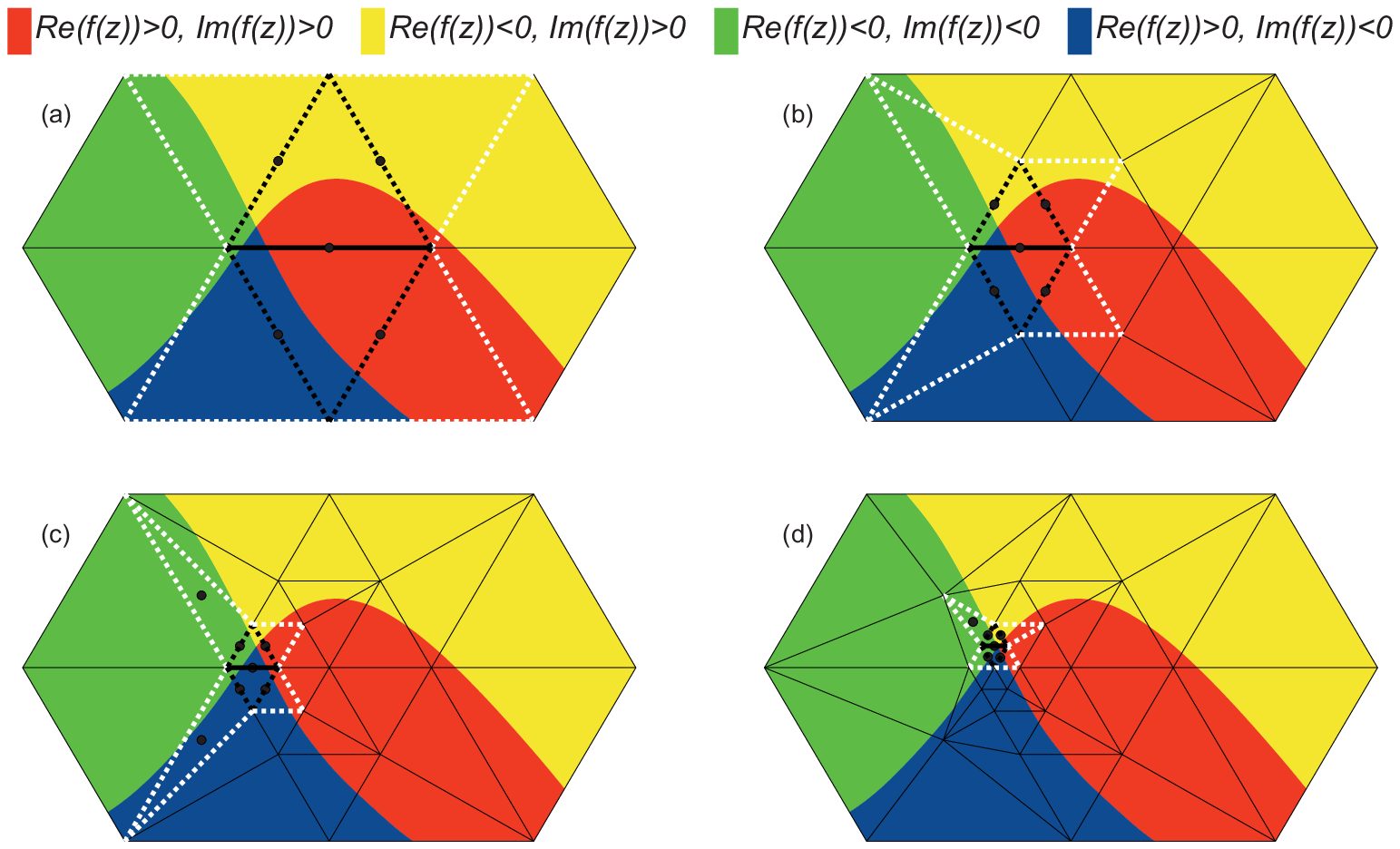}}
\caption{A simple example of the mesh refinement process. Figures (a)-(d) represent four consecutive iterations. The candidate edges are marked by thick black lines. The black dotted lines represent the boundary of the candidate regions and the white dotted line represents a boundary of the extra zone.} \label{fig:E2}
\end{figure}

In order to illustrate the main idea of the proposed approach, a simple example of the process is presented in Figure \ref{fig:E2} (again, a phase portrait of the function is placed in the background). In the first step, new extra nodes are added to the mesh in the centers of the edges in the candidate regions. Then, using Delaunay triangulation, a new mesh is obtained. Next, the function values are evaluated at these new points and the new configuration is analyzed exactly as in the preliminary estimation - new candidate regions are determined for a locally denser mesh. Obviously, the area of the new candidate region is smaller, which improves the accuracy of the result. The process may be repeated any number of times, until a fixed accuracy $\delta$ is reached. 

In subsequent repetitions, the refinement of the mesh can lead to ill-conditioned geometry ("skinny triangles"). To avoid this problem an additional zone surrounding the region should be considered (white dotted line in Figure \ref{fig:E2}). If the triangle in the extra zone is "skinny" (e.g. the ratio of the longest triangle edge to the shortest edge is greater than $3$), a new extra node is added in its center - see Figures \ref{fig:E2} (c) and \ref{fig:E2} (d).

In this stage of the algorithm, the refinement can be performed only for roots or poles (if there is no need to find all the characteristic points). However, it is possible (and quite efficient) to start the refinement process without verification of the candidate regions in the preliminary estimation. The verification can be performed after the refinement, and the results may be more accurate (e.g., two different roots could be verified as a double root in the preliminary estimation, but they may be separated in the refinement process).

\subsection{Effectiveness and Limitations of the Algorithm}\label{sec:limits}
An application of the regular mesh has a very clear guarantee of correctness -  
if the discretization of the function is proper, then none of the roots/poles can be missed. The proper discretization means that for all edges the phase change does not exceed three quadrants. Hence, the algorithm can be applied for any analytic function and in any domain, if the initial mesh step $\Delta r$ is sufficiently small. 

However, if the self-adaptive process is involved this condition is sufficient, but not necessary (for the initial mesh discretization). In practice, $\Delta r$ can be even greater than the distance between the roots/poles of the function (as shown in section \ref{sec:num}). Unfortunately, just as for the other established methods (e.g., based on interpolatins \cite{Austin2014} or discrete techniques \cite{Kowalczyk15,Wan11,Meylan2003}), there is no clear recipe for the a priori estimation of the initial sampling, for an arbitrary function. In practice, it can be chosen by a user experimentally, via sequential iterations. Initial verification can be performed using the idea of DCAP for the whole boundary of the region $\Omega$. However, this still does not guarantee that all the roots are found (for instance if the result of DCAP is equal to zero, then the region may be free of roots and poles or it may contain an equal number of roots and poles). 

To reduce the risk of missing roots/poles, CAP can be extended to higher moments $m$ \cite{Lamparillo75}:
\begin{equation}
\label{eqn:cap_himo} \frac{1}{2 \pi i}\oint\limits_{C} \frac{f'(z)}{f(z)} z^m dz =\sum_{k\in\{roots\}}\left( z^{(k)}\right)^m - \sum_{k\in\{poles\}}\left( z^{(k)}\right)^m.
\end{equation}

The moment $m=1$ eliminates the problem for a single root-pole pair \cite{Zieniutycz83} and each higher moment can further reduce the risk. Such an approach can be especially useful if an analytical expression of the function is known.

\section{Numerical Tests}\label{sec:num}
The algorithm was implemented in the MATLAB environment, and all the tests were performed using an Intel(R) Core i7-2600K CPU 3.40-GHz, 16-GB RAM computer.

\subsection{Complex Modes}\label{sec:cmpx}
As the first example, a complex wave propagation problem in a
circular waveguide of radius $b$, coaxially loaded with a dielectric cylinder of radius $a$ is considered \cite{Mrozowski97,Kowalczyk17AP}. To ensure continuity of the fields at the boundary of the dielectric and metal, the following determinant function must
be equal to zero:
\begin{equation}
\label{eqn:cmpx} f(z)= \left|
\begin{array}{cccccc}
-J_1 & 0 & J_2 & Y_2& 0  &0   \\
0   &J_1& 0    & 0   &-J_2&-Y_2\\
-\frac{z m J_1}{a\kappa_1^2}& -\frac{i\eta_0 J'_1}{\kappa_1}& \frac{z m J_2}{a\kappa_2^2}& \frac{z m Y_2}{a\kappa_2^2}&\frac{i\eta_0 J'_2}{\kappa_2}&\frac{i\eta_0 Y'_2}{\kappa_2}\\
-\frac{i\varepsilon_r J'_1}{\kappa_1\eta_0}& -\frac{z m J_1}{a\kappa_1^2}&
\frac{i J'_2}{\kappa_2\eta_0}& \frac{i Y'_2}{\kappa_2\eta_0}& \frac{z m J_2}{a\kappa_2^2}&
\frac{z m Y_2}{a\kappa_2^2}\\
0& 0& J_3& Y_3& 0&0\\
0& 0& \frac{z m J_3}{b\kappa_2}& \frac{ z m Y_3}{b\kappa_2}& i\eta_0J'_3&
i\eta_0 Y'_3\\
\end{array} \right|,
\end{equation}
where $z$ represents a normalized propagation coefficient. $J_1=J_m(k_0\kappa_1 a)$, $Y_1=Y_m(k_0\kappa_1 a)$, $J_2=J_m(k_0\kappa_2 a)$, $Y_2=Y_m(k_0\kappa_2 a)$, $J_3=J_m(k_0\kappa_2 b)$ and $Y_3=Y_m(k_0\kappa_2 b)$  are Bessel and Neumann functions (primes denote derivatives). The coefficients are $\kappa_1=\sqrt{z^2+\varepsilon_r}$, $\kappa_2=\sqrt{z^2+1}$, $k_0=2\pi f / c$ and $\eta_0=120\pi$ $\Omega$. The tests are performed for parameters 
$a = 6.35$ mm, $b = 10$ mm, $\varepsilon_r=10$, $m=1$ and $f=5$ GHz. 

\begin{figure}
\centerline{\includegraphics[width=1.0\columnwidth]{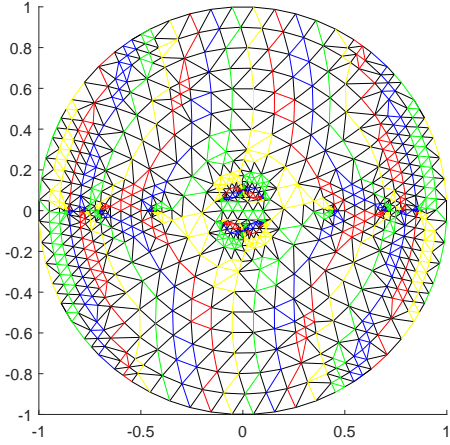}}
\caption{Self-adaptive meshes obtained for problem (\ref{eqn:cmpx}); initial mesh $N=271$ and final mesh $N=6167$ for $\delta=1e-15$.} \label{fig:cmpx}
\end{figure}

In order to compare the efficiency of the proposed algorithm with the other established methods, the considered region  is a unite disk $\Omega=\{\bar{z}\in\mathbb{C}:|\bar{z}|<1\}$ and the scaling factor $10$ is applied $z=10\bar{z}$. 
The initial mesh, evenly covering region $\Omega$ with $N=271$ nodes, is sufficient to find all roots and poles of the function in the considered region. Such discretization corresponds to mesh resolution $\Delta r=0.15$ and, obviously, any higher resolution leads to the same results - twelve single roots:\\
$\bar{z}^{(1)}=-0.096642302459942 - 0.062923397455697i$,\\
$\bar{z}^{(2)}=-0.096642302459942 + 0.062923397455697i$,\\
$\bar{z}^{(3)}=0.096642302459942 - 0.062923397455697i$,\\
$\bar{z}^{(4)}=0.096642302459942 + 0.062923397455696i$,\\
$\bar{z}^{(5)}=-0.444429043110023 + 0.000000000000000i$,\\
$\bar{z}^{(6)}=0.444429043110023 - 0.000000000000000i$,\\
$\bar{z}^{(7)}=-0.703772250217811 + 0.000000000000000i$,\\
$\bar{z}^{(8)}=0.703772250217811 - 0.000000000000000i$,\\
$\bar{z}^{(9)}=-0.775021522202022 + 0.000000000000000i$,\\
$\bar{z}^{(10)}=0.775021522202023 - 0.000000000000000i$,\\
$\bar{z}^{(11)}=-0.856115203911565 + 0.000000000000000i$,\\
$\bar{z}^{(12)}=0.856115203911564 - 0.000000000000000i$\\
and two second order poles:\\
$\bar{z}^{(13)}=0.000000000000000 + 0.100000000000000i$,\\
$\bar{z}^{(14)}=0.000000000000000 - 0.100000000000000i$.\\
The parameters and results of the analysis for various accuracy $\delta$ are collected  in Table \ref{tab:cmpx} and in Figure \ref{fig:cmpx}. 

\begin{table}%
\caption{Analysis of problem (\ref{eqn:cmpx}) - parameters and results of the proposed algorithm \label{tab:cmpx}}
\begin{center}
\begin{tabular}{|c|c|c|c|}
\hline
accuracy & CPU time & no. of nodes & no. of iterations\\
\hline
$\delta = {1e-3}$ & $0.46$ s  & $1603$ & $11$  \\
\hline
$\delta = {1e-6}$ & $0.72$ s  & $2759$ & $20$  \\
\hline
$\delta = {1e-9}$ & $1.05$ s  & $3867$ & $30$  \\
\hline
$\delta = {1e-12}$ & $1.37$ s  & $5013$ & $40$  \\
\hline
$\delta = {1e-15}$ & $1.78$ s  & $6167$ & $51$  \\
\hline
\end{tabular}
\end{center}
\end{table}%

Similar discrete algorithm \cite{Kowalczyk15} requires $721$ samples (corresponds to $\Delta r=0.06$) for the initial mesh and this number increases to $2252$ nodes for the verification and for the accuracy improvement to level $1e-3$. Moreover, the computational time of algorithm \cite{Kowalczyk15} is two orders of magnitude longer and equals $113$ s. This large discrepancy in the algorithms arises from fundamental differences in the data processing. The older algorithm \cite{Kowalczyk15} is much more complex. It requires approximation of the real and imaginary parts of the function separately. Next, the curves representing the zero of the real and imaginary parts of the function are constructed, and then all crossings of these curves are estimated. Finally, DCAP (with extra nodes) is applied over the circles surrounding the crossings. Moreover, if DCAP is applied over the artificial circle of radius $\Delta r$ surrounding the candidate point, this value $\Delta r$ must be sufficiently small to separate all the roots, which requires much denser initial discretization. In the approach presented in this paper, the processing is significantly simpler, which results in faster calculations and lower memory requirements.

Additionally, the effectiveness of the method based on rational approximation presented in \cite{Austin2014}, based on the "ratdisk" algorithm \cite{Gonnet2011}, is tested. The parameters and results of the analysis for various sampling $N$ and orders of numerator $m$ and denominator polynomials $n$ are collected in Table \ref{tab:cmpx_rd}. If the number of samples $N$ is sufficiently high, then all the roots can be found. However, despite the regularization \cite{Gonnet2011} used in the analysis, there are numerous spurious roots and poles among the proper results. Moreover, the accuracy of the results cannot be simply controlled - for higher number of samples $N$ and higher orders of $m$ and $n$ the accuracy increases, but for each root/pole this value can be different (maximum $Err_{max}$ and minimum $Err_{min}$ absolute errors are collected in Table \ref{tab:cmpx_rd}). So, in practice, these results can be used as an efficient preliminary estimation and extra post processing may be required.

\begin{table}%
\caption{Analysis of problem (\ref{eqn:cmpx})  - parameters and results of rational approximation \cite{Austin2014}  \label{tab:cmpx_rd}}
\begin{center}
\begin{tabular}{|p{0.6cm}|p{1.8cm}|p{1.8cm}|p{1.8cm}|}
\hline
            & $N=50$             & $N=500$           & $N=5000$\\
\hline
$m=25$      & $Err_{max}=1e-1$   & $Err_{max}=1e-1$  & $Err_{max}=1e-1$   \\
$n=25$      & $Err_{min}=2e-4$   & $Err_{min}=7e-4$  & $Err_{min}=7e-4$   \\
            & $t_{CPU}=0.09s$    & $t_{CPU}=0.18s$   & $t_{CPU}=0.81s$     \\
            & $2$ missing roots  & $2$ missing roots & $2$ missing roots   \\
            & $1$ spurious pole  &                   &                     \\            
\hline
$m=250$     &     $m+n>N$        & $Err_{max}=4e-3$  & $Err_{max}=2e-3$   \\
$n=250$     &  					 & $Err_{min}=8e-9$  & $Err_{min}=6e-9$   \\
            & 					 & $t_{CPU}=0.25s$   & $t_{CPU}=1.32s$     \\
            & 					 & $7$ spurious roots& $4$ spurious roots  \\
            & 					 & $7$ spurious poles& $6$ spurious poles  \\            
\hline
$m=2500$    &     $m+n>N$        &     $m+n>N$      & $Err_{max}=2e-3$   \\
$n=2500$    &  					 & 					& $Err_{min}=5e-9$   \\
            & 					 & 					& $t_{CPU}=23.08s$     \\
            & 					 & 					& $3$ spurious roots  \\
            & 					 & 					& $5$ spurious poles  \\            
\hline
\end{tabular}
\end{center}
\end{table}%

\subsection{Lossy Multilayered Waveguide}\label{sec:multi}
As the second example, a multilayered guiding structure is considered \cite{Anemogiannis92,ZOUROS2018}. That structures are widely used in microwave applications and their analysis boils down to satisfying specific boundary conditions, which requires zero of the following determinant function:
\begin{equation}
\label{eqn:multi} 
f(z)= \left|
\begin{array}{cc}
1  & -\cos(k_0\kappa_1 d_1)-\gamma_C \sin(k_0 \kappa_1 d_1)/\kappa_1   \\
i\gamma_S  & -i \kappa_1\sin(k_0\kappa_1 d_1)+i\gamma_C \cos(k_0 \kappa_1 d_1)   \\
\end{array} \right|,
\end{equation}
where $z$ represents a normalized propagation coefficient, $k_0=2\pi /\lambda_0$, $\kappa_1=\sqrt{n_1^2-z^2}$, $\gamma_S=\sqrt{z^2-n_S^2}$ and $\gamma_C=\sqrt{z^2-n_C^2}$. The typical set of material parameters is $n_1=1.5835$, $n_S=0.065-4i$ and $n_C=1$. The analysis is performed for thickness $d_1=1.81$ $\mu$m and frequency corresponds to wavelength $\lambda_0=0.6328$ $\mu$m.

\begin{figure}
\centerline{\includegraphics[width=1.0\columnwidth]{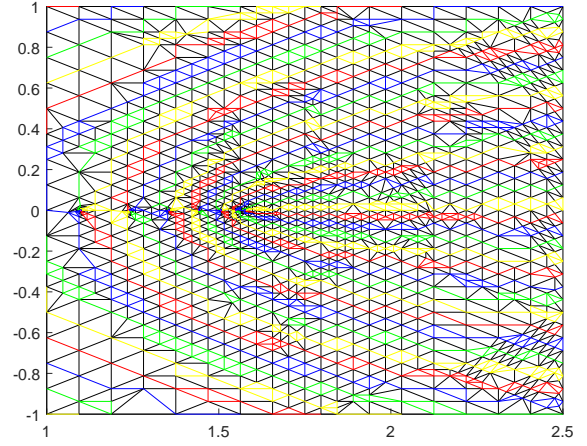}}
\caption{Self-adaptive meshes obtained for problem (\ref{eqn:multi}), initial mesh $N=27$, final mesh $N=3322$ for $\delta=1e-15$.} \label{fig:multi}
\end{figure}

The assumed region is the same as the one proposed in \cite{ZOUROS2018} $\Omega=\{z\in\mathbb{C}:1<\textrm{Re}(z)<2.5 \wedge -1<\textrm{Im}(z)<1\}$. The results and parameters of the analysis for various accuracy $\delta$ are shown in Table \ref{tab:multi} and in Figure \ref{fig:multi}. 

In this case, the initial mesh, evenly covering region $\Omega$ with $N=27$ nodes (which corresponds to $\Delta r=0.5$), is sufficient to find all roots of the function (\ref{eqn:multi}) in this region  - seven single roots:\\
$z^{(1)}=1.574863045752781 - 0.000002974623699i$,\\
$z^{(2)}=1.548692243882210 - 0.000012101013332i$,\\
$z^{(3)}=1.504169866404311 - 0.000028029436583i$,\\
$z^{(4)}=1.439795544245059 - 0.000052001665381i$,\\
$z^{(5)}=1.353140429182476 - 0.000086139194522i$,\\
$z^{(6)}=1.240454471356097 - 0.000133822149870i$,\\
$z^{(7)}=1.096752543407689 - 0.000197146879192i$.

Again, the efficiency is compared with discrete algorithm \cite{Kowalczyk15}, which requires $10927$ ($\Delta r=0.02$) samples for the initial mesh and this number increases to $11067$ nodes for the verification and for the accuracy improvement to level $1e-3$. Also, the computational time of the algorithm \cite{Kowalczyk15} is about two orders of magnitude longer and equals $62.33$  s. The most recently published algorithm \cite{ZOUROS2018} requires even more function calls (in this case it is $156803$) which results in significantly longer analysis. The same applies to algorithm \cite{Anemogiannis92} where the similar huge number of the function samples is required.

\begin{table}%
\caption{Analysis of problem (\ref{eqn:multi}) - parameters and results of the proposed algorithm \label{tab:multi} }
\begin{center}
\begin{tabular}{|c|c|c|c|}
\hline
accuracy & CPU time & no. of nodes & no. of iterations\\
\hline
$\delta = {1e-3}$ & $0.33$  s  & $1623$ & $10$  \\
\hline
$\delta = {1e-6}$ & $0.44$  s  & $2066$ & $21$  \\
\hline
$\delta = {1e-9}$ & $0.58$  s  & $2472$ & $31$  \\
\hline
$\delta = {1e-12}$ & $0.71$  s  & $2900$ & $41$  \\
\hline
$\delta = {1e-15}$ & $0.87$  s  & $3322$ & $51$  \\
\hline
\end{tabular}
\end{center}
\end{table}%

\subsection{Graphene Transmission Line}\label{sec:gtl}
As the last example, a simple graphene transmission line is considered. The guide consists of a thin graphene layer deposited and a silicone substrate \cite{Gomez2013,Kowalczyk17AP}. In this case the normalized propagation coefficient $z$, for TM modes, can be found from the following equation
\begin{eqnarray}
\label{eqn:gtl} 
&& f(z)=
\frac{\varepsilon_{r1}}{\eta_0 \sqrt{\varepsilon_{r1}+z^2}}+
\frac{\varepsilon_{r2}}{\eta_0 \sqrt{\varepsilon_{r2}+z^2}}\nonumber\\
&&\qquad +\left[ 
\sigma_{lo}-z^2k_{0}^{2}(\alpha_{sd}+\beta_{sd})
\right], 
\end{eqnarray}
where $k_0=2\pi f/c$ and $\eta_0$ is a wave impedance of the vacuum. The graphene parameters depend on the frequency as follows
\begin{equation}
\sigma_{lo}=\frac{-iq_e^2k_B T}{\pi \hbar^2 (2\pi f-i\tau^{-1})} 
\ln\left[ 
2\left(1+\cosh \left( \frac{\mu_c}{k_B T} \right)    \right)
\right],
\end{equation}
\begin{equation}
\alpha_{sd}=\frac{-3 v_F^2 \sigma_{lo} }{4 (2\pi f-i\tau^{-1})^2}, \qquad \beta_{sd}=\frac{\alpha_{sd}}{3},
\end{equation}
where $q_e$ is electron charge, $k_B$ is Boltzmann's constant, $T=300$ K, $\tau=0.135$ ps,  
$\mu_c=0.05q_e$, $v_F=10^6$ m/s. The tests are performed for frequency $f=1$ THz, $\varepsilon_{r1}=1$ and  $\varepsilon_{r2}=11.9$.

The region $\Omega=\{z\in\mathbb{C}:-100<\textrm{Re}(z)<400 \wedge -100<\textrm{Im}(z)<400\}$ is considered. Due to four Riemann sheets of the function (\ref{eqn:gtl}) their pointwise product is analyzed \cite{Kowalczyk17}, in order to avoid separate investigation of each sheet.

\begin{figure}
\centerline{\includegraphics[width=1.0\columnwidth]{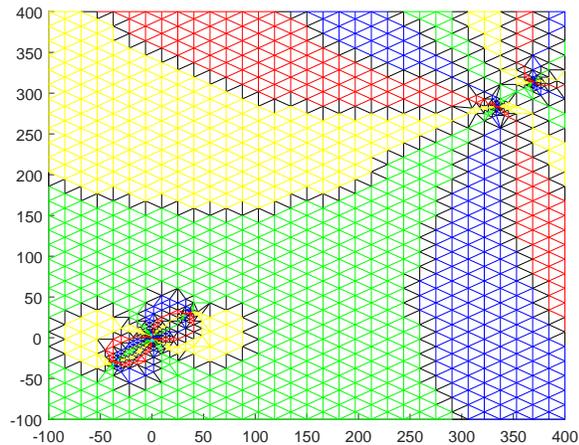}}
\caption{Self-adaptive meshes obtained for problem (\ref{eqn:gtl}), initial mesh $N=973$, final mesh $N=4983$ for $\delta = {1e-12}$.} \label{fig:gtl}
\end{figure}

The results and parameters of the analysis for various accuracy $\delta$ are presented in Table \ref{tab:gtl} and in Figure \ref{fig:gtl}. $N=973$ function samples ($\Delta r=18$) is sufficient to determine all roots and poles of the function in $\Omega$ - eight single roots:\\
$z^{(1)}= -32.1019622516073 - 27.4308619360125i$,\\
$z^{(2)}= 32.1019622516073 + 27.4308619360128i$,\\
$z^{(3)}=-38.1777253144799 - 32.5295210455987i $,\\
$z^{(4)}= 38.1777253144797 - 32.5295210455987i $,\\
$z^{(5)}= 332.7448889298402 + 282.2430799544401i$,\\
$z^{(6)}= 336.2202873389791 + 285.1910910139915i$,\\
$z^{(7)}= 368.4394672155518 + 312.5220780593669i $,\\
$z^{(8)}= 371.0075708341529 + 314.7004076766967i$,\\
and two second order poles:\\ 
$z^{(9)}=0.000000000000184 - 3.449637662132114i$,\\
$z^{(10)}=-0.000000000000158 + 3.449637662131965i$.

\begin{table}%
\caption{Analysis of problem (\ref{eqn:gtl}) - parameters and results \label{tab:gtl} }
\begin{center}
\begin{tabular}{|c|c|c|c|}
\hline
accuracy & CPU time & no. of nodes & no. of iterations\\
\hline
$\delta = {1e-3}$ & $0.39$ s  & $2342$ & $16$  \\
\hline
$\delta = {1e-6}$ & $0.54$ s  & $3121$ & $26$  \\
\hline
$\delta = {1e-9}$ & $0.75$ s  & $4084$ & $36$  \\
\hline
$\delta = {1e-12}$ & $0.99$ s  & $4983$ & $46$  \\
\hline
\end{tabular}
\end{center}
\end{table}%

For comparison, the algorithm \cite{Kowalczyk15} requires $32689$ (corresponds to $\Delta r=3$) samples for the initial mesh and this number increases to $33256$ nodes for the verification and for the accuracy improvement to level $1e-3$. In this case, the computational time is even three orders of magnitude longer and equals $487$ s.

\section{Conclusions}
A new algorithm for complex roots and poles finding is presented. A wide class of complex functions can be analyzed in any arbitrarily shaped search region. The effectiveness of the proposed technique is supported by several numerical tests. Moreover, the efficiency of the presented method is confirmed by comparing the analysis parameters to those obtained from alternative recently published techniques. It is shown that the proposed algorithm can be up to three orders of magnitude faster and requires significantly smaller number of function evaluations.

\appendices
\section{Source code}
The source code for the GRPF (Global complex Roots and Poles Finding algorithm based on phase analysis), can be found at: (if the paper is accepted, the code will be available at https://github.com/), and it is licensed under the MIT License.

\begin{IEEEbiographynophoto}{Piotr Kowalczyk}
was born in Wejherowo, Poland, in 1977. He received the M.S. degree in applied physics and mathematics and Ph.D. degree in electrical engineering from the Gdansk University of Technology, Gdansk, Poland, in 2001 and 2008, respectively.
He is currently with Microwave and Antenna Engineering, Technical University,
Gdansk, Poland. His research is focused on scattering and propagation of electromagnetic wave problems, algorithms and numerical methods.
\end{IEEEbiographynophoto}

\end{document}